\documentclass[11pt,a4paper]{article}
\usepackage{mathrsfs}
\usepackage{amssymb}
\usepackage{amsfonts}
\usepackage{color,xcolor}
\usepackage{graphicx}
\usepackage{manfnt}
\RequirePackage[colorlinks,citecolor=blue,urlcolor=blue]{hyperref}
\newtheorem{theorem}{Theorem}[section]
\newtheorem{corollary}{Corollary}[section]
\newtheorem{lemma}{Lemma}[section]

\newtheorem{definition}{Definition}[section]

\def\<{{\langle}}\def\>{{\rangle}}
\def\[{{\Big[}}\def\]{{\Big]}}\def\({{\Big(}}\def\){{\Big)}}

\def\={&\!\!=\!\!&}
\def\cA{{\mathcal A}}

\def\cS{{\mathcal S}}

\def\mN{{\mathbb N}}

\def\mZ{{\mathbb Z}}

\def\geq{\geqslant}\def\leq{\leqslant}
\def\inf{\infty}
\def\div{\mathord{{\rm div}}}

\begin{document}
\title{\bf Existence and uniqueness of $W^{1,r}_{loc}$-solutions for stochastic transport equations\footnote{This work was partly supported by the NSFC grants 11501577, 11531006, 11771449, 11771123 and PAPD of Jiangsu Higher Education Institutions.}}

\author{Jinlong Wei$^a$, Jinqiao Duan$^b$, Hongjun Gao$^c$ and Guangying Lv$^d$}
\date{{$^a$ School of Statistics and Mathematics, Zhongnan University of}\\
{Economics and Law, Wuhan, Hubei 430073, China}\\
{weijinlong@zuel.edu.cn}
\\ {$^b$ Department of Applied Mathematics}\\
{Illinois Institute of Technology, Chicago, IL 60616, USA}\\
{duan@iit.edu} \\  {$^c$ Institute of Mathematics, School of Mathematical Science} \\ {Nanjing Normal University, Nanjing 210023, China}\\ {gaohj@njnu.edu.cn}
\\ {$^d$ Institute of Contemporary Mathematics, Henan University}\\
{Kaifeng, Henan 475001, China} \\
{gylvmaths@henu.edu.cn}}

 \maketitle
\noindent{\hrulefill}

\vskip1mm\noindent{\bf Abstract} We investigate a stochastic transport
equation driven by a multiplicative noise. For $L^q(0,T;W^{1,p}({\mathbb R}^d;{\mathbb R}^d))$ drift coefficient and $W^{1,r}({\mathbb R}^d)$ initial data, we obtain the existence and uniqueness of stochastic strong solutions (in $W^{1,r}_{loc}({\mathbb R}^d))$.
In particular, when $r=\infty$, we establish a Lipschitz estimate for solutions and this question is opened by Fedrizzi and Flandoli in case of  $L^q(0,T;L^p({\mathbb R}^d;{\mathbb R}^d))$ drift coefficient. Moreover, opposite to the
deterministic case where $L^q(0,T;W^{1,p}({\mathbb R}^d;{\mathbb R}^d))$ drift coefficient and $W^{1,p}({\mathbb R}^d)$ initial data may induce non-existence for strong
solutions (in $W^{1,p}_{loc}({\mathbb R}^d)$), we prove that a multiplicative stochastic perturbation of
Brownian type is enough to render the equation well-posed. It
is an interesting example of a deterministic PDE that
becomes well-posed under the influence of a multiplicative Brownian
type noise. We extend the existing results \cite{FF2,FGP1} partially.



 \vskip0mm\noindent{\hrulefill}
\section{Introduction}\label{sec1}\setcounter{equation}{0}
In recent years the theory of stochastic partial differential
equations has had an intensive development and many important
contributions have been obtained \cite{DaP}. One of the important
branch in this field is to touch the effect of noise on the
solutions. Opposite to the deterministic case, we find that a
multiplicative stochastic perturbation of Brownian type will
discover new phenomena of stochastic nature. For example, noise can
make the weak solutions unique in linear transport equation
\cite{AF,FGP1,NO}, can make solutions more regular \cite{BFG}, can
prevent singularities in linear transport equations \cite{FF2,FGP3}, can
prevent infinite stretching of the passive field in a stochastic
vector advection equation \cite{FMN}, can prevent collapse of
Vlasov-Poisson point charges \cite{DFV}. Besides a multiplicative
stochastic perturbation, while an additive stochastic perturbation
or random initial value may induce new properties for solutions of
differential equations. For example, noise can regularize the
solutions for ordinary differential
equations \cite{Att,FL,FF1,FF3,Fl,FGP2,KR,Zhang}.

In this paper, we will discuss this topic and study the effect of noise on the solutions to a linear transport
equation. To be precise, given $T>0$, let us consider the following Cauchy problem
\begin{eqnarray}\label{1.1}
\left\{
  \begin{array}{ll}
\partial_tu(t,x)+b(t,x)\cdot\nabla u(t,x)
+\sum_{i=1}^d\partial_{x_i}u(t,x)\circ\dot{B}_i(t)=0, \ (\omega,t,x)\in \Omega\times (0,T)\times {\mathbb R}^d, \\
u(t,x)|_{t=0}=u_0(x), \  x\in{\mathbb R}^d,
  \end{array}
\right.
\end{eqnarray}
where $B(t)=(B_1(t), B_2(t), _{\cdots}, B_d(t))$ is a
$d$-dimensional standard Brownian motion defined on a stochastic
basis ($\Omega, \mathcal{F},{\mathbb P},(\mathcal{F}_{t})_{t\geq 0}$), the stochastic
integration with a notation $\circ$ is interpreted in Stratonovich
sense, $b: [0,T]\times{\mathbb R}^d\rightarrow{\mathbb R}^d, \ u_0:
{\mathbb R}^d\rightarrow{\mathbb R}$ are measurable functions which are in class of $L^1(0,T;L^1_{loc}({\mathbb R}^d;{\mathbb R}^d))$ and $L^1(0,T;L^1_{loc}({\mathbb R}^d))$ respectively. We are going to investigate stochastic strong solutions. Here the stochastic strong is defined as the following
\begin{definition} \label{def1.1} Let $r\in [1,\infty]$, $\div b\in L^1(0,T;L^{r^\prime}_{loc}({\mathbb R}^d))$ ($1/r+1/r^\prime=1$). Suppose that $u\in L^\infty(\Omega\times[0,T];L^r({\mathbb R}^d))$ is a random field.
We call $u$ a stochastic weak solution of (\ref{1.1}) if for every
$\varphi\in\mathcal{C}_0^\infty({\mathbb R}^d)$, $\int_{{\mathbb R}^d}\varphi(x)u(t,x)dx$
has a continuous modification which is an $\mathcal{F}_t$-semimartingale and
for every $t\in [0,T]$,
\begin{eqnarray}\label{1.2}
\int_{{\mathbb R}^d}\varphi(x)u(t,x)dx&=&\int_{{\mathbb R}^d}\varphi(x)u_0(x)dx+
\int^t_0\int_{{\mathbb R}^d}\div(b(s,x)\varphi(x))u(s,x)dxds\cr\cr&&
+\sum_{i=1}^d\int^t_0\circ dB_i(s)\int_{{\mathbb R}^d}\partial_{x_i}\varphi(x)u(s,x)dx,  \quad
{\mathbb P}-a.s..
\end{eqnarray}
Moreover, if
\begin{eqnarray}\label{1.3}
{\mathbb P}\Big(\ |\nabla u|\in L^\infty(0,T;L^r_{loc}({\mathbb R}^d))\Big)=1,
\end{eqnarray}
then $u$ is called a stochastic strong solution.
\end{definition}

Our first main result is
\begin{theorem} \label{the1.1} \textbf{(Existence and uniqueness)} Let $q\in (2,\infty]$, $p\in [2,\infty)$, such that
\begin{eqnarray}\label{1.4}
\frac{2}{q}+\frac{d}{p}<1, \ \ b\in L^q(0,T;W^{1,p}({\mathbb R}^d;{\mathbb R}^d)), \ \ \div b\in L^1(0,T;L^\infty({\mathbb R}^d)).
\end{eqnarray}
We assume that $r\in [1,\infty]$ and $u_0\in L^r({\mathbb R}^d), |\nabla u_0|\in L^r_{loc}({\mathbb R}^d)$. Then the following statements hold:

(i) there exists a unique stochastic strong solution to the Cauchy problem (\ref{1.1}). Moreover, the unique
stochastic strong solution can be represented by $u(t,x)=u_0(X^{-1}(t,x))$. Here $\{X(t,x)\}$ is the unique
strong solution of (\ref{2.1}) with $s=0$;

(ii) when $r=\infty$, if the initial $u_0$ also yields that $|\nabla u_0|\in L^\infty({\mathbb R}^d)$, then the unique strong solution belongs to $\cap_{a\geq 1}L^a(\Omega;W^{1,\infty}_{loc}([0,T]\times {\mathbb R}^d))$.
\end{theorem}

\vskip2mm\noindent
\textbf{Remark 1.1.} (i) In \cite{FF2} Fedrizzi and Flandoli also discussed the problem (\ref{1.1}), under the assumptions that
\begin{eqnarray*}
b\in L^q(0,T;L^p({\mathbb R}^d;{\mathbb R}^d)), \ q,p\in [2,\infty), \ \frac{2}{q}+\frac{d}{p}<1, \ u_0\in \cap_{r\geq 1}W^{1,r}({\mathbb R}^d),
\end{eqnarray*}
they conclude the existence and uniqueness of $\cap_{r\geq 1}W^{1,r}_{loc}$ solutions, i.e. for every $t\in [0,T]$,
\begin{eqnarray}\label{1.5}
{\mathbb P}\Big(u(t)\in \cap_{r\geq 1}W^{1,r}_{loc}({\mathbb R}^d)\Big)=1.
\end{eqnarray}
Then an interesting questions is posed:

$\bullet$ How to extend the persistence property of solutions for a fixed $r\in [1,\infty]$, i.e. if $u_0\in W^{1,r}({\mathbb R}^d)$, does we have
\begin{eqnarray*}
{\mathbb P}\Big(u(t)\in W^{1,r}_{loc}({\mathbb R}^d)\Big)=1 \ ?
\end{eqnarray*}
(posed by Fedrizzi and Flandoli in \cite{FF2}: "After the result of Theorem 1 (conclusion (\ref{1.5})), it remains open the question whether the solution is Lipschitz
continuous (or more) when $u_0\in W^{1,\infty}({\mathbb R}^d)$ (or more)").

And it is the major source of our present paper. From Theorem \ref{the1.1} we gave a positive answer for above question. However, here we suppose in addition that $|\nabla b|\in L^q(0,T;L^p({\mathbb R}^d))$ ($\div b\in L^1(0,T;L^\infty({\mathbb R}^d))$ is not necessary when $r=\infty$).

(ii) When proving Theorem \ref{the1.1}, the key component is to show that (\ref{1.3}) holds true. We give it two proofs. The first one is based on the continuity of the stochastic field $\nabla_xX^{-1}(t,x)$ in $(t,x)$ (or see Lemma \ref{lem2.1}). However, it seems impossible to carry over this approach to the case of stochastic partial differential equations (such as stochastic transport equations) driven by non-Gaussian L\'{e}vy processes directly since now the stochastic field is not continuous in $t$. The second proof is based upon a moment estimate on $\|\nabla_xX^{-1}(t,x)\|_{L^\infty((0,T)\times B_R)}$ (for every $R>0$).
It turns out that this device works equally well for stochastic partial differential equations driven by  non-Gaussian L\'{e}vy processes and we will use the device to deal with stochastic transport equations driven by $\alpha$-stable L\'{e}vy processes in future.

(iii) For the drift coefficient in critical $L^q(0,T;L^p({\mathbb R}^d;{\mathbb R}^d))$ space, i.e.
\begin{eqnarray*}
b\in L^q(0,T;L^p({\mathbb R}^d;{\mathbb R}^d)), \ q,p\in [2,\infty), \ \frac{2}{q}+\frac{d}{p}\leq 1, \ or \ p=\infty, \ q=2,
\end{eqnarray*}
and $u_0\in \cap_{r\geq 1}W^{1,r}({\mathbb R}^d)$, Beck, Flandoli, Gubinelli and Maurelli \cite{BFG} proved that (\ref{1.5}) holds true as well. However, in this critical case, we do not know whether or not  (\ref{1.3}) is valid yet.

(iv) For the drift coefficient is bounded, i.e. $q=p=\infty$ and $u_0\in \mathcal{C}^1_b$, Mohammed, Nilssen and Proske \cite{MNP} showed that (\ref{1.1}) exists a unique bounded weak solution. Moreover, the authors proved that for every $t>0$ and every $p\in (1,\infty)$, $u(t)\in W^{1,p}({\mathbb R}^d,w)$ (a weighted Sobolev space). Different from this result, under the assumption (\ref{1.4}), we show that if the initial data is Lipschitz continuous, then the unique solution is also Lipschitz continuous (but only locally).

In view of Theorem \ref{the1.1}, we conclude the following comparison principle.
\begin{corollary}\label{cor1.1} \textbf{(Comparison Principle)} Let $p,q,r,d$ and $b$ be described in Theorem \ref{the1.1}
and $u_{0,1},u_{0,2}\in L^r({\mathbb R}^d)$, $|\nabla u_{0,1}|,|\nabla u_{0,2}|\in L^r_{loc}({\mathbb R}^d)$ . Assume that $u_1$ and $u_2$
are two stochastic strong solutions of (\ref{1.1}), with initial
values $u_{0,1}$ and $u_{0,2}$. If $u_{0,1}\leq u_{0,2}$, then with
probability 1, $u_1\leq u_2$. In particular, if the initial value is
nonnegative, then with probability 1, the unique stochastic strong
solution is nonnegative as well.
\end{corollary}

In applications to nonlinear equations, especially for Euler and Navier-Stokes equations, $b$ should maintain the same Sobolev regularity as the solutions. From this point, we have the following persistence property of solutions.

\begin{corollary}\label{cor1.2} \textbf{(Persistence property)} Let $p>d$ such that $u_0\in L^p({\mathbb R}^d)$, $|\nabla u_0|\in L^p_{loc}({\mathbb R}^d)$, $b\in L^\infty(0,T;W^{1,p}({\mathbb R}^d;{\mathbb R}^d))$ and $\div b\in L^1(0,T;L^\infty({\mathbb R}^d))$.
Then there exists a unique stochastic strong solution to the Cauchy problem (\ref{1.1}) (i.e. $r$ is replaced by $p$ in Theorems \ref{the1.1}). Moreover, the unique stochastic strong solution can be represented by $u(t,x)=u_0(X^{-1}(t,x))$.
\end{corollary}

The scope of work is not limited to describe a well-posed result, we also present a counter example for non-existence of such solutions in the deterministic case. Precisely speaking, when the noise vanishes, we prove that now the strong solution will not be Sobolev differentiable.
\begin{theorem} \label{the1.2} \textbf{(Non-existence)} Let $p,b$ and $u_0$ be stated in Corollary \ref{cor1.2}. Consider the Cauchy problem
\begin{eqnarray}\label{1.6}
\left\{
  \begin{array}{ll}
\partial_tu(t,x)+b(t,x)\cdot\nabla u(t,x)=0, \ (t,x)\in (0,T)\times {\mathbb R}^d, \\
u(t,x)|_{t=0}=u_0(x), \  x\in{\mathbb R}^d.
  \end{array}
\right.
\end{eqnarray}
Then there exists a unique weak solution. Moreover, the unique weak
solution can be represented by $u(t,x)=u_0(X^{-1}(t,x))$. Here
$X(t,x)$ is the unique DiPerna-Lions flow dominated by the ODE
\begin{eqnarray}\label{1.7}
dX(t)=b(t,X(t))dt, \ \ t\in(0,T], \ \ X(t)|_{t=0}=x.
\end{eqnarray}
However, if $d\geq 2$, one can choose proper functions $u_0$ and $b$ such
that $u(t,x)=u_0(X^{-1}(t,x))$ does not lie in
$L^\infty(0,T;W^{1,p}_{loc}({\mathbb R}^d)).$ Here $u$ is said to be a weak
solution of (\ref{1.6}), if it lies in $L^\infty(0,T;L^p({\mathbb R}^d))$
and for every $\varphi\in\mathcal{C}_0^\infty({\mathbb R}^d)$, every $t\in [0,T]$,
\begin{eqnarray}\label{1.8}
\int_{{\mathbb R}^d}\varphi(x)u(t,x)dx=\int_{{\mathbb R}^d}\varphi(x)u_0(x)dx+
\int^t_0\int_{{\mathbb R}^d}\div(b(s,x)\varphi(x))u(s,x)dxds.
\end{eqnarray}
If $|\nabla u|\in L^\infty(0,T;L^p_{loc}({\mathbb R}^d))$ in addition, we call $u$ a strong solution of (\ref{1.6}).
 \end{theorem}

\vskip2mm\noindent
\textbf{Remark 1.2.} The equation in (\ref{1.6}) is hyperbolic, even though for smooth initial value, solutions spontaneously develop discontinuities. However, here we prove that the emergence of shocks is prevented when one adds a multiplicative Brownian type noise which preserves the hyperbolic structure of the equation.

The rest of the paper is structured as follows. In Sections 2-3, the
proof of Theorems \ref{the1.1} and \ref{the1.2} are given. Section 2
is concerned with the existence and uniqueness of the stochastic
strong solutions and a counter example for non-existence is given in
Section 3.

\vskip2mm\noindent
\textbf{Notations}  The letter $C$ will mean a positive constant, whose values may change in different places. For a parameter or a function $\varrho$, $C(\varrho)$ means the constant is only dependent on $\varrho$. $\mN$ is the set of natural numbers and $\mZ$ denotes the set of integral numbers. For every $R>0$, $B_R:=\{x\in{\mathbb R}^d:|x|<R\}$. Almost surely can be abbreviated to $a.s.$.

\section{Proof of Theorem \ref{the1.1}}\label{sec3}
\setcounter{equation}{0}

In this section, we shall give the proof of Theorem \ref{the1.1}. Firstly, we give a definition and a useful lemma.

Given $s\in[0,T]$ and $x\in {\mathbb R}^d$, consider the stochastic differential equation (SDE for short) in ${\mathbb R}^d$:
\begin{eqnarray}\label{2.1}
dX(s,t)=b(t,X(s,t))dt+dB(t), \ \ t\in(s,T], \ \ X(s,t)|_{t=s}=x.
\end{eqnarray}

\begin{definition} [\cite{Kun2}, P114] \label{def2.1} A stochastic homeomorphism flow (respect. of class $\mathcal{C}^{1,\beta}$ with $\beta\in (0,1)$) on
$(\Omega, \mathcal{F},{\mathbb P}, (\mathcal{F}_t)_{0\leq t\leq T})$ associated to (\ref{2.1}) is a
map $(s,t,x,\omega) \rightarrow X(s,t,x)(\omega)$, defined for
$0\leq s \leq t \leq T, \ x\in {\mathbb R}^d, \ \omega \in \Omega$ with
values in ${\mathbb R}^d$, such that

(i) given any $s\in [0,T],\  x \in {\mathbb R}^d$, the process
$\{X(s,\cdot,x)\}= \{X(s,t,x), \ t\in [s,T]\}$ is a continuous
$\{\mathcal{F}_{s,t}\}_{s\leq t\leq T}$-adapted solution of (\ref{2.1});

(ii) ${\mathbb P}-a.s.$, for all  $0\leq s \leq t \leq T$, the functions
$X(s,t,x), \ X^{-1}(s,t,x)$ are continuous in $(s,t,x)$;

(iii) ${\mathbb P}-a.s., \   X(s,t,x)=X(r,t,X(s,r,x))$  for all
$0\leq s\leq r \leq t \leq T$, $x\in {\mathbb R}^d$ and $X(s,s,x)=x$.
\end{definition}

For further use, we review a useful lemma.
\begin{lemma} [\cite{FGP1}] \label{lem2.1} Let us assume that there exists $\alpha \in (0,1)$ such
that $b\in L^{\infty}(0,T; \mathcal{C}^\alpha_b({\mathbb R}^d;{\mathbb R}^d))$. Then for
every $s\in [0,T ], \  x\in {\mathbb R}^d$, the stochastic equation
(\ref{2.1}) has a unique continuous adapted solution $\{X(s,t,x), \
t\in [s,T ], \ \omega \in\Omega\}$, which forms a ${\mathcal
C}^{1,\alpha^\prime}$ ($\alpha^\prime<\alpha$) stochastic flow $X(s,t)$ of diffeomorphisms. Moreover, if we let $(b^n)\subset L^{\infty}(0,T;\mathcal{C}^{\alpha}_b({\mathbb R}^d ;{\mathbb R}^d))$ be a sequence of vector
fields and $X^n$ be the corresponding stochastic flows. If $b^n
\rightarrow b$ in $L^{\infty}(0,T; \mathcal{C}^{\alpha^\prime}_b({\mathbb R}^d ;{\mathbb R}^d))$ for some $\alpha^\prime> 0$, then for any $a\geq 1$,
\begin{eqnarray}\label{2.2}
\lim _{n\rightarrow\infty}\sup_{x\in {\mathbb R}^d} \sup_{0\leq s\leq
T}{\mathbb E}[\sup_{r\in[s,T]}|X^n(s,r,x)-X(s,r,x)|^a]=0,
\end{eqnarray}
\begin{eqnarray}\label{2.3}
\sup_{n}\sup_{x\in {\mathbb R}^d}\sup_{0\leq s\leq
T}{\mathbb E}[\sup_{r\in[s,T]}\|DX^n(s,r,x)\|^a]<\infty,
\end{eqnarray}
\begin{eqnarray}\label{2.4}
\lim_{n\rightarrow \infty}\sup_{x\in {\mathbb R}^d} \sup_{0\leq s\leq
T}{\mathbb E}[\sup_{r\in[s,T]}\|DX^n(s,r,x)-DX(s,r,x)\|^a]=0,
\end{eqnarray}
(recall that $\| \cdot \|$ is the Hilbert-Schmidt norm on $d \times
d$ matrices).
\end{lemma}

\vskip2mm\noindent \textbf{Remark 2.1.} When the coefficients $b^n$ are smooth
enough, using the Liouville theorem, we have Euler's
identity:
\begin{eqnarray*}
det(\nabla_xX^n(t,x))=\exp(\int^t_0\div b^n(r,X^n(r,x)))dr).
\end{eqnarray*}
In view of (\ref{2.2}) and (\ref{2.4}), if one suppose in addition
that $b^n\rightarrow b$ in $L^1(0,T; W^{1,a}_{loc}({\mathbb R}^d;{\mathbb R}^d))$
for some $a\geq 1$, up to choosing a subsequence, one derives
\begin{eqnarray}\label{2.5}
det(\nabla_xX(t,x))=\exp(\int^t_0\div b(r,X(r,x)))dr),
\end{eqnarray}
here for simplicity, we have used $X(t,x)$ to stand for $X(0,t,x)$.

\vskip2mm\noindent \textbf{Proof of Theorem \ref{the1.1}.}  (i) Initially, we prove the uniqueness and noticing that the equation is linear, we need to show $u\equiv 0$ a.s. if the initial value vanishes. For $r=\infty$, the uniqueness has been proved by Attanasio and Flandoli \cite[Theorem 11]{AF}. Therefore, it suffices to show $r<\infty$. Let $\varrho_\varepsilon$ be a regularizing kernel i.e.
\begin{eqnarray*}
\varrho_\varepsilon =\frac{1}{\varepsilon^d} \varrho(\frac{\cdot}{\varepsilon}) \ \ with \ \ 0\leq \varrho \in \mathcal{C}^\infty_0({\mathbb R}^d) , \ \ support(\varrho)\subset B_1, \ \ \int_{{\mathbb R}^d}\varrho(x)dx=1.
\end{eqnarray*}
We define $u_\varepsilon=u\ast \varrho_\varepsilon$, then $u_\varepsilon$ yields that
\begin{eqnarray*}
\partial_tu_\varepsilon(t,x)+b(t,x)\cdot\nabla u_\varepsilon(t,x)
+\sum_{i=1}^d\partial_{x_i}u_\varepsilon(t,x)\circ\dot{B}_i(t)=I_\varepsilon,
\end{eqnarray*}
with
\begin{eqnarray*}
I_\varepsilon=b(t,x)\cdot\nabla u_\varepsilon(t,x)-(b\cdot\nabla u)_\varepsilon(t,x).
\end{eqnarray*}
With the help of assumption (\ref{1.4}),
\begin{eqnarray*}
I_\varepsilon\rightarrow 0 \ \ in \ \ L^q(0,T;L^r_{loc}({\mathbb R}^d)), \quad {\mathbb P}-a.s..
\end{eqnarray*}
By approximation arguments (see \cite[Corollary II.1]{DL} but with some minor modifications), for every $M>0$ and $r<\infty$ one ends up with
\begin{eqnarray*}
\frac{d}{dt}\int_{{\mathbb R}^d}|u(t,x)\wedge M|^rdx\leq C \int_{{\mathbb R}^d}|u(t,x)\wedge M|^rdx, \quad {\mathbb P}-a.s..
\end{eqnarray*}
From this one proves the uniqueness.

Secondly, we show that $u(t,x)=u_0(X^{-1}(t,x))$ is a stochastic weak solution of (\ref{1.1}). Here $\{X(t,x)\}$ is the unique
strong solution of (\ref{2.1}) with $s=0$. We divide the proof into two cases: $q=\infty, d<p<\infty$ and $p\in [2,\infty), q\in (2,\infty)$, $2/q+d/p<1$. We begin our proof for the first case.
\vskip2mm\par
\textbf{$\bullet$ Case 1:} $q=\infty, d<p<\infty$.
\vskip2mm\par
Since $p>d$, $L^\infty(0,T;W^{1,p}({\mathbb R}^d;{\mathbb R}^d)) \subseteq
L^\infty(0,T;\mathcal{C}_b^\alpha({\mathbb R}^d;{\mathbb R}^d))$ with $\alpha=1-d/p$. By Lemma \ref{lem2.1}, the stochastic differential equation
(\ref{2.1}) with $s=0$ has a unique continuous adapted solution $\{X(t,x), \ t\in [0,T ], \ \omega \in\Omega\}$, which forms a $\mathcal{C}^{1,\alpha^\prime}$ ($\alpha^\prime<\alpha$) stochastic flow $X(t,x)$ of diffeomorphisms. If one defines $u(t,x)=u_0(X^{-1}(t,x))$ and uses the Kunita-It\^{o}-Wentzel formula (see \cite[Theorem 8.3]{Kun1} or \cite[Lemma 2.1]{CO}), for every $\varphi\in\mathcal{C}_0^\infty({\mathbb R}^d)$, $\int_{{\mathbb R}^d}\varphi(x)u(t,x)dx$ meets (\ref{1.2}). Thus $\int_{{\mathbb R}^d}\varphi(x)u(t,x)dx$ has a continuous modification which is an $\mathcal{F}_t$-semimartingale. To
complete the proof, we need to show $u\in
L^\infty(\Omega\times(0,T);L^r({\mathbb R}^d))$ for $r\in [1,\infty]$. Clearly when $r=\infty$, it is true. It remains to show $r\in [1,\infty)$.

With the help of (\ref{2.5}), then
\begin{eqnarray}\label{2.6}
\int_{{\mathbb R}^d}|u_0(X^{-1}(t,x))|^rdx&=&\int_{{\mathbb R}^d}|u_0(x)|^rdet(\nabla_xX(t,x))dx
\cr\cr&=&\int_{{\mathbb R}^d}|u_0(x)|^r\exp(\int^t_0\div b(\tau,X(\tau,x))d\tau)dx\cr\cr&\leq&
\exp(\|\div b\|_{L^1(0,T;L^\infty({\mathbb R}^d))})\int_{{\mathbb R}^d}|u_0(x)|^rdx.
\end{eqnarray}

\textbf{$\bullet$ Case 2:} $p\in [2,\infty), q\in (2,\infty)$, $\frac{2}{q}+\frac{d}{p}<1$.
\vskip2mm\par
Recall that
SDE (\ref{2.1}) is equivalent to the following SDE (see
\cite{FF2,FF3,FGP1}):
\begin{eqnarray}\label{2.7}
d Y(t)=\lambda U(t,\gamma^{-1}(t,Y(t)))dt+ [I+\nabla
U(t,\gamma^{-1}(t,Y(t)))] dB(t), \ t\in(0,T], \ Y(t)|_{t=0}=y,
\end{eqnarray}
where $\gamma(t,x)=x+U(t,x)$, $\gamma^{-1}(t,x)$ is its inverse of
the mapping $x\mapsto\gamma(t,x)$, $U$ is given by
\begin{eqnarray}\label{2.8}
\left\{
\begin{array}{ll}
\partial_{t}U(t,x) +\frac{1}{2}\Delta U(t,x)+b(t,x)\cdot \nabla U(t,x)=
\lambda U(t,x)-b(t,x), \ (t,x)\in (0,T)\times {\mathbb R}^d, \\
U(T,x)=0, \  x\in{\mathbb R}^d.
  \end{array}
\right.
\end{eqnarray}
In addition the solutions (\ref{2.7}) and (\ref{2.1}) has the
relationship $X(t)=\gamma^{-1}(t)\circ Y(t)$. With the aid of \cite[Theorem 3.3, Lemma 3.4]{FF3} (or see \cite[Theorem 1.2]{Kry}), there is a unique $U\in L^q(0,T;W^{2,p}({\mathbb R}^d;{\mathbb R}^d))\cap W^{1,q}(0,T;L^p({\mathbb R}^d;{\mathbb R}^d))$ solving the backward PDE (\ref{2.8}). Moreover, there is a finite constant $N$ such that
\begin{eqnarray}\label{2.9}
\|\partial_tU\|_{L^q(0,T;L^p({\mathbb R}^d))} + \|U\|_{L^q(0,T;W^{2,p}({\mathbb R}^d))}\leq N\|b\|_{L^q(0,T;L^p({\mathbb R}^d))}.
\end{eqnarray}
From (\ref{2.9}), by the Morrey inequality (see \cite{Eva}, P282), the Sobolev embedding inequality (see \cite{Eva}, P306), (\ref{1.4}) and the maximum principle then
\begin{eqnarray}\label{2.10}
\|\nabla U\|_{L^\infty((0,T)\times {\mathbb R}^d)} \leq N\|b\|_{L^q(0,T;L^p({\mathbb R}^d))}, \ \  \|\nabla U\|_{L^\infty((0,T)\times {\mathbb R}^d)} \rightarrow 0, \ as \ \lambda\rightarrow \infty.
\end{eqnarray}

For $1\leq i\leq d$, if one differentiates $x_i$ in equation (\ref{2.8}) and denotes $U_1=\partial_{x_i}U$ then
\begin{eqnarray*}
\left\{
\begin{array}{ll}
\partial_{t}U_1(t,x) +\frac{1}{2}\Delta U_1(t,x)+b(t,x)\cdot \nabla U_1(t,x)=
\lambda U_1(t,x)-\partial_{x_i}b(t,x)-\partial_{x_i}b(t,x)\cdot \nabla U(t,x), \\
U_1(T,x)=0.
  \end{array}
\right.
\end{eqnarray*}
Thus $U_1\in L^q(0,T;W^{2,p}({\mathbb R}^d;{\mathbb R}^d))\cap W^{1,q}(0,T;L^p({\mathbb R}^d;{\mathbb R}^d))$. Since $1\leq i\leq d$, one concludes that $\nabla U \in L^q(0,T;W^{2,p}({\mathbb R}^d;{\mathbb R}^{d\times d}))\cap W^{1,q}(0,T;L^p({\mathbb R}^d;{\mathbb R}^{d\times d}))$ and by (\ref{2.10}),
\begin{eqnarray}\label{2.11}
\|\nabla^2 U\|_{L^\infty((0,T)\times {\mathbb R}^d)} \leq N\|\nabla b\|_{L^q(0,T;L^p({\mathbb R}^d))}(1 +N\|b\|_{L^q(0,T;L^p({\mathbb R}^d))}).
\end{eqnarray}

Combining (\ref{2.10}) and (\ref{2.11}), for a given and big enough real number $\lambda$, from \cite[Lemma 3.5]{FF3}, for every $t\in[0,T]$, $\gamma(t,x)=x+U(t,x)$ forms a non-singular diffeomorphism of $\mathcal{C}^2$. Moreover, $\gamma$ and $\gamma^{-1}$ have a bounded first and second spatial derives and $\nabla^2 \gamma$ is globally H\"{o}lder continuous. Therefore, (\ref{2.7}) exists a unique strong solution which forms a $\mathcal{C}^{1,\beta}$ ($0<\beta<1$) stochastic flow of diffeomorphisms. The remains is the same as discussed in the case of $q=\infty, d<p<\infty$. We achieve the proof.

Thirdly, we show that (\ref{1.3}) holds. Noticing that both $q=\infty, d<p<\infty$ and $p\in [2,\infty), q\in (2,\infty)$, $2/q+d/p<1$, the stochastic differential equation (\ref{2.1}) beginning from $s=0$ has a unique continuous adapted solution $\{X(t,x), \ t\in [0,T ], \ \omega \in\Omega\}$, which forms a $\mathcal{C}^{1,\alpha^\prime}$ ($0<\alpha^\prime<\alpha$) stochastic flow $X(t,x)$ of diffeomorphisms.
We have the following chain rule
\begin{eqnarray}\label{2.12}
\nabla_x(u_0(X^{-1}(t,x)))=\nabla_xu_0(X^{-1}(t,x))
\nabla_xX^{-1}(t,x),
\end{eqnarray}
which implies that: for every $R>0$ and $r\in[1,\infty)$,
\begin{eqnarray}\label{2.13}
\int_{B_R}|\nabla_x(u_0(X^{-1}(t,x)))|^rdx=
\int_{B_R}|\nabla_xu_0(X^{-1}(t,x))|^r\|\nabla_xX^{-1}(t,x)\|^rdx,
\end{eqnarray}
where $\| \cdot \|$ is the Hilbert-Schmidt norm on $d \times d$
matrices. We will show that the right hand in (\ref{2.13}) is finite almost surely. To reach this aim, let us give it two proofs.

 \vskip2mm\par
\textbf{$\bullet$ Case 1:} $q=\infty, d<p<\infty$.
\vskip2mm\par
\textbf{The first proof.} We only recall the idea of the proof, see \cite{FGP3} for details. Since now $b\in L^\infty(0,T;\mathcal{C}_b^\alpha({\mathbb R}^d;{\mathbb R}^d))$, $X(t,x)$ and $\nabla_xX^{-1}(t,x)$ are continuous in $(t,x)$ almost surely. From (\ref{2.12}), (\ref{2.13}) and Euler's identity (\ref{2.5}), one fulfills our conclusion.

\textbf{The second proof.}  The proof is based on estimating ${\mathbb E}\sup_{0\leq t\leq T, x\in B_R}\|\nabla_xX^{-1}(t,x)\|$ and the calculations can be divided into four steps.
\vskip2mm\par
\textbf{Step 1.} To simplify the problem.
\vskip2mm\par

Since the backward flow satisfies the same SDE of the forward flow with a drift coefficient of opposite, for every $t\in [0,T]$, to calculate ${\mathbb E}\sup_{0\leq t\leq T, x\in
B_R}\|\nabla_xX^{-1}(t,x)\|$ it is sufficient to estimate
${\mathbb E}\sup_{0\leq t\leq T, x\in B_R}\|\nabla_xX(t,x)\|$. Recall that
SDE (\ref{2.1}) is equivalent to (\ref{2.7}), and  $X(t)=\gamma^{-1}(t)\circ Y(t)$. Using \cite[Lemma 3.5]{FF3}, then one gets that $\|\nabla\gamma^{-1}(t)\|\leq 2$, hence if one can manipulate ${\mathbb E}\sup_{0\leq t\leq T, y\in B_R}\|\nabla_yY(t,y)\|$, then we accomplish our proof. On the other hand, by scaling transforming: $y=2Ry_1$ first, and shift transforming $y_1=x_1+(1/2,_{\cdots},1/2) (\in {\mathbb R}^d)$, we need to show ${\mathbb E}\sup_{0\leq t\leq T, x_1\in [0,1]^d}\|\nabla Y(t,x_1)\|$. For notation no confusion, in the following calculation, one also write $x$ instead of $x_1$.
\vskip2mm\par
\textbf{Step 2.} Space H\"{o}lder estimates for $Y(t)$.
\vskip2mm\par

Let $Y(t,x)$ and $Y(t,y)$ be the unique strong solution of
(\ref{2.7}) with initial data $x$ and $y$ respectively. If one sets
$Y_t(x,y)=Y(t,x)-Y(t,y)$ and $\tilde{b}(t,y)=\lambda
U(t,\gamma^{-1}(t,y))$, $\sigma(t,y)=I+\nabla
U(t,\gamma^{-1}(t,y))$, then
\begin{eqnarray*}
\left\{\begin{array}{ll} d
Y_t(x,y)=[\tilde{b}(t,Y(t,x))-\tilde{b}(t,Y(t,y))]dt
+[\sigma(t,Y(t,x))-\sigma(t,Y(t,x))]d B(t), \ t\in(0,T), \\
Y_t(x,y)|_{t=0}=x-y.
\end{array}
\right.
\end{eqnarray*}
Using the It\^{o} formula, for $m\geq 2$, we have
\begin{eqnarray}\label{2.14}
&&d |Y_t(x,y)|^m\cr\cr&=&m|Y_t(x,y)|^{m-2}\langle Y_t(x,y),
\tilde{b}(t,Y(t,x))-\tilde{b}(t,Y(t,y))\rangle dt\cr\cr&&+
\frac{1}{2}m(m-1)|Y_t(x,y)|^{m-2}tr([\sigma(t,Y(t,x))-\sigma(t,Y(t,x))][\sigma(t,Y(t,x))-\sigma(t,Y(t,x))]^\top)dt
\cr\cr&&+m|Y_t(x,y)|^{m-2}\langle Y_t(x,y),
[\sigma(t,Y(t,x))-\sigma(t,Y(t,x))]d B(t)\rangle.
\end{eqnarray}
Observing that $U$ is the unique solution of (\ref{2.8}) and $b\in
L^\infty(0,T;W^{1,p}({\mathbb R}^d))$, so $U\in
L^\infty(0,T;\mathcal{C}^{2+\alpha}_b({\mathbb R}^d))$ $(\alpha=1-d/p)$, and then
$\tilde{b} \in L^\infty(0,T;Lip({\mathbb R}^d;{\mathbb R}^d))$, $\sigma \in L^\infty(0,T;Lip({\mathbb R}^d;{\mathbb R}^{d\times d}))$. We obtain from
(\ref{2.14}) that
\begin{eqnarray*}
d |Y_t(x,y)|^m&\leq&C(m)|Y_t(x,y)|^m
dt\cr\cr&&+m|Y_t(x,y)|^{m-2}\langle Y_t(x,y),
[\sigma(t,Y(t,x))-\sigma(t,Y(t,x))]d B(t)\rangle,
\end{eqnarray*}
i.e.
\begin{eqnarray}\label{2.15}
|Y_t(x,y)|^m&\leq&|x-y|^m+C(m)\int_0^t|Y_s(x,y)|^m
ds\cr&&+m\int_0^t|Y_s(x,y)|^{m-2}\langle Y_s(x,y),
[\sigma(s,Y(s,x))-\sigma(s,Y(s,x))]d B(s)\rangle.
\end{eqnarray}
Therefore
\begin{eqnarray}\label{2.16}
{\mathbb E}|Y_t(x,y)|^m\leq|x-y|^m+C(m)\int_0^t{\mathbb E}|Y_s(x,y)|^mds.
\end{eqnarray}
From (\ref{2.16}), if one uses the Gr\"{o}nwall inequality, it
yields that
\begin{eqnarray}\label{2.17}
\sup_{0\leq t\leq T}{\mathbb E}|Y_t(x,y)|^m\leq C(m,T)|x-y|^m.
\end{eqnarray}

On the other hand, by virtue of the BDG inequality, from
(\ref{2.15}), one concludes that
\begin{eqnarray*}
{\mathbb E}\sup_{0\leq s\leq
t}|Y_s(x,y)|^m&\leq&|x-y|^m+C(m)\int_0^t{\mathbb E}\sup_{0\leq r\leq
s}|Y_r(x,y)|^m dr\cr&&+C(m){\mathbb E}\[\int_0^t|Y_s(x,y)|^{2m}ds\]^{\frac{1}{2}}.
\end{eqnarray*}
Since (\ref{2.17}) holds for every $m\geq 2$, by (\ref{2.17}) and
Minkowshi's inequality, then
\begin{eqnarray*}
{\mathbb E}\sup_{0\leq s\leq t}|Y_s(x,y)|^m\leq
C(m,T)|x-y|^m+C(m)\int_0^t{\mathbb E}\sup_{0\leq r\leq s}|Y_r(x,y)|^m dr,
\end{eqnarray*}
which suggests that
\begin{eqnarray*}
{\mathbb E}\sup_{0\leq t\leq T}|Y_t(x,y)|^m\leq C(m,T)|x-y|^m,
\end{eqnarray*}
if one uses the Gr\"{o}nwall inequality again. From this, one
finishes at
\begin{eqnarray}\label{2.18}
{\mathbb E}\sup_{0\leq t\leq T}|X_t(x,y)|^m\leq C(m,T)|x-y|^m.
\end{eqnarray}

\vskip2mm\par
\textbf{Step 3.}  H\"{o}lder estimate for $\|\nabla_xY(t,x)\|$.
\vskip2mm\par

Let us set $\nabla_xY(t,x)$ by $\xi_t(x)$, from  (\ref{2.7}) then for every $t\in (0,T]$, $\xi_t(x)$ yields that
\begin{eqnarray*}
d \xi_t(x)=\lambda \nabla U(t,X(t,x))\nabla\gamma^{-1}(t,Y(t))\xi_t(x)dt+\nabla^2
U(t,X(t,x))\nabla\gamma^{-1}(t,Y(t))\xi_t(x)d B(t),
\end{eqnarray*}
and $\xi_t(x)|_{t=0}=I$.

By It\^{o}'s formula, for every $m\geq 2$, we have
\begin{eqnarray}\label{2.19}
&&d \|\xi_t(x)\|^m\cr\cr&=&m\|\xi_t(x)\|^{m-2}\langle\xi_t(x), \lambda
\nabla U(t,X(t,x))\nabla\gamma^{-1}(t,Y(t))\xi_t(x)\rangle dt\cr\cr&&+
\frac{1}{2}m(m-1)\|\xi_t(x)\|^{m-2}tr([\nabla^2
U(t,X(t,x))\nabla\gamma^{-1}(t,Y(t))\xi_t(x)]
\cr\cr&&\cdot[\nabla^2U(t,X(t,x))\nabla\gamma^{-1}(t,Y(t))\xi_t(x)]^\top)dt
\cr\cr&&+m\|\xi_t(x)\|^{m-2}\langle\xi_t(x),
\nabla^2 U(t,X(t,x))\nabla\gamma^{-1}(t,Y(t))\xi_t(x)d B(t)\rangle .
\end{eqnarray}

Observing that $U\in L^\infty(0,T;\mathcal{C}^{2+\alpha}_b({\mathbb R}^d))$
$(\alpha=1-d/p)$. From (\ref{2.19}), one fulfills that
\begin{eqnarray*}
d \|\xi_t(x)\|^m\leq
C\|\xi_t(x)\|^m+m\|\xi_t(x)\|^{m-2}\langle\xi_t(x), \nabla^2
U(t,X(t,x))\nabla\gamma^{-1}(t,Y(t))\xi_t(x)d B(t)\rangle,
\end{eqnarray*}
which suggests
\begin{eqnarray}\label{2.20}
{\mathbb E}\|\xi_t(x)\|^m\leq C+ C\int^t_0{\mathbb E}\|\xi_s(x)\|^mds.
\end{eqnarray}

From (\ref{2.20}), one applies the Gr\"{o}nwall inequality to get
\begin{eqnarray}\label{2.21}
\sup_{0\leq t\leq T, x\in{\mathbb R}^d}{\mathbb E}\|\xi_t(x)\|^m\leq C.
\end{eqnarray}
The calculations from (\ref{2.17}) to (\ref{2.18}) uses here again, suggests that
\begin{eqnarray}\label{2.22}
\sup_{x\in{\mathbb R}^d}{\mathbb E}\sup_{0\leq t\leq T}\|\xi_t(x)\|^m\leq C.
\end{eqnarray}

If one set $\xi_t(x,y)=\xi_t(x)-\xi_t(y)$, by an analogue manipulation of (\ref{2.14}), it reaches at
\begin{eqnarray*}
d \|\xi_t(x,y)\|^m&=&m\|\xi_t(x,y)\|^{m-2}\langle\xi_t(x,y), A_1\rangle
dt+m\|\xi_t(x,y)\|^{m-2}\langle\xi_t(x,y), A_2d B(t)\rangle
\cr\cr&&+ \frac{1}{2}m(m-1)\|\xi_t(x,y)\|^{m-2}tr(A_3A_3^\top)dt,
\end{eqnarray*}
where
\begin{eqnarray*}
&&A_1=\lambda
\nabla U(t,X(t,x))\nabla\gamma^{-1}(t,Y(t,x))\xi_t(x)-\lambda \nabla U(t,X(t,y))\nabla\gamma^{-1}(t,Y(t,y))\xi_t(y),
\cr\cr
&&A_2=\nabla^2
U(t,X(t,x))\nabla\gamma^{-1}(t,Y(t,x))\xi_t(x)-\nabla^2 U(t,X(t,y))\nabla\gamma^{-1}(t,Y(t,y))\xi_t(y),
\cr\cr
&&A_3=\nabla^2 U(t,X(t,x))\nabla\gamma^{-1}(t,Y(t,x))\xi_t(x)-\nabla^2
U(t,X(t,y))\nabla\gamma^{-1}(t,Y(t,y))\xi_t(y),
\end{eqnarray*}
which implies
\begin{eqnarray*}
&&d \|\xi_t(x,y)\|^m
\cr\cr&\leq&m\|\xi_t(x,y)\|^{m-1}\|A_1\| dt+
\frac{1}{2}m(m-1)\|\xi_t(x,y)\|^{m-2}\| A_3\|^2dt \cr\cr&&+
m\|\xi_t(x,y)\|^{m-2}\langle\xi_t(x,y), A_2d B(t)\rangle
\cr\cr &\leq&C(m)\|\xi_t(x,y)\|^{m-1}\Big(\|\xi_t(x,y)\|+
\|\xi_t(x)\|[|X(t,x)-X(t,y)|+|X(t,x)-X(t,y)|^\alpha]\Big)
dt\cr\cr&&+ m\|\xi_t(x,y)\|^{m-2}\langle\xi_t(x,y), A_2d B(t)\rangle.
\end{eqnarray*}
Therefore
\begin{eqnarray}\label{2.23}
{\mathbb E}\|\xi_t(x,y)\|^m &\leq&C(m)\int_0^t{\mathbb E}\|\xi_s(x,y)\|^mds\cr\cr&&+
C(m){\mathbb E}\int_0^t\|\xi_s(x)\|^m[|X(s,x)-X(s,y)|^m+|X(s,x)-X(s,y)|^{\alpha
m}]ds \cr\cr&\leq& C(m)\int_0^t{\mathbb E}\|\xi_s(x,y)\|^mds \cr\cr&&+
C(m)\int_0^t\Big({\mathbb E}\|\xi_s(x)\|^{2m}\Big)^{\frac{1}{2}}\Big(
{\mathbb E}[|X(s,x)-X(s,y)|^{2m}\Big)^{\frac{1}{2}}ds \cr\cr&&+
C(m)\int_0^t\Big({\mathbb E}\|\xi_s(x)\|^{2m}\Big)^{\frac{1}{2}}\Big(
{\mathbb E}|X(s,x)-X(s,y)|^{2\alpha m}\Big)^{\frac{1}{2}}ds .
\end{eqnarray}
Combining (\ref{2.18}) and (\ref{2.21}), (\ref{2.23}) induces that
\begin{eqnarray}\label{2.24}
{\mathbb E}\|\xi_t(x,y)\|^m \leq C(m)\int_0^t{\mathbb E}\|\xi_s(x,y)\|^mds+
C(m,T)[|x-y|^m+|x-y|^{\alpha m}].
\end{eqnarray}
Thus
\begin{eqnarray}\label{2.25}
\sup_{0\leq t\leq T}{\mathbb E}\|\xi_t(x,y)\|^m \leq
C(m,T)[|x-y|^m+|x-y|^{\alpha m}].
\end{eqnarray}

Similar manipulations of (\ref{2.17})-(\ref{2.18}) applies again,
one ends up with
\begin{eqnarray*}
&&{\mathbb E}\sup_{0\leq s\leq t}\|\xi_s(x,y)\|^m\cr\cr&\leq&
C(m)\int_0^t{\mathbb E}\sup_{0\leq r\leq s}\|\xi_s(x,y)\|^mds \cr\cr&&+
C(m){\mathbb E}\int_0^t\|\xi_s(x)\|^m[|X(s,x)-X(s,y)|^m+|X(s,x)-X(s,y)|^{\alpha
m}]ds \cr\cr&& + C(m){\mathbb E}\Big[\int_0^t\|\xi_s(x,y)\|^{2m}
ds\Big]^{1/2} \cr\cr&&+ C(m){\mathbb E}\Big[\int_0^t \|\xi_s(x)\|^{2m}|X(s,x)-X(s,y)|^{2\alpha m} ds\Big]^{\frac{1}{2}}
\cr\cr&\leq&
C(m)\int_0^t{\mathbb E}\sup_{0\leq r\leq s}\|\xi_s(x,y)\|^mds \cr\cr&&+
C(m)\int_0^t\Big[{\mathbb E}\|\xi_s(x)\|^{2m}\]^{\frac{1}{2}}
\Big[{\mathbb E}|X(s,x)-X(s,y)|^{2m}+{\mathbb E}|X(s,x)-X(s,y)|^{2\alpha m}]^{\frac{1}{2}} ds
\cr\cr&& + C(m)\Big[\int_0^t{\mathbb E}\|\xi_s(x,y)\|^{2m}
ds\Big]^{1/2} \cr\cr&&+ C(m)\Big[\int_0^t \Big({\mathbb E}\|\xi_s(x)\|^{4m}\Big)^{\frac{1}{2}} \Big({\mathbb E}|X(s,x)-X(s,y)|^{4\alpha m}\Big)^{\frac{1}{2}} ds\Big]^{\frac{1}{2}}
\cr\cr&\leq&
C(m)\int_0^t{\mathbb E}\sup_{0\leq r\leq s}\|\xi_s(x,y)\|^mds+
C(m,T)[|x-y|^m+|x-y|^{\alpha m}].
\end{eqnarray*}
The Gr\"{o}nwall inequality uses again, we gain
\begin{eqnarray}\label{2.26}
&&{\mathbb E}\sup_{0\leq t\leq T}\|\xi_t(x,y)\|^m\cr\cr&\leq&
C(m,T)[|x-y|^m+|x-y|^{\alpha m}]\leq C(m,T)|x-y|^{\alpha m}, \ \forall \ x,y\in [0,1]^d.
\end{eqnarray}

\vskip2mm\par
\textbf{Step 4.}  ${\mathbb E}\sup_{0\leq t\leq T, x\in
[0,1]^d}\|\nabla_xY(t,x)\|<\infty.$
\vskip2mm\par

To arrive our purpose, we introduce a sequence of sets: $\cS_n=\{z\in
\mZ^d\ | \ z2^{-n}\in [0,1]^d\}, \ n\in \mN$. For an arbitrary
$e=(e^1,_{\cdots},e^d)\in \mZ^d$ such that $\|e\|_\infty=\max_{1\leq
i\leq d}|e^i|=1$, and every $z,z+e\in\cS_n$, we define
$\xi_z^{n,e}(t)=|\xi_t((z+e)2^{-n})-\xi_t(z2^{-n})|$. Then by
(\ref{2.26}), for every $m\geq 2$,
\begin{eqnarray*}
{\mathbb E}\sup_{0\leq t\leq T}|\xi_z^{n,e}(t)|^m\leq C(m,T)2^{-n \alpha m}.
\end{eqnarray*}
For any $\tau>0$ and $K>0$, one sets a number of events
$\cA_{z,\tau}^{n,e}=\{\omega\in\Omega \ | \ \sup_{0\leq t\leq
T}\xi_z^{n,e}(t)\geq K\tau^n\}$ ($z,z+e\in\cS_n$), it yields that
\begin{eqnarray*}
{\mathbb P}(\cA_{z,\tau}^{n,e})\leq \frac{{\mathbb E}\sup_{0\leq t\leq
T}|\xi_z^{n,e}(t)|^m}{K^m\tau^{mn}}\leq \frac{C(m,T)2^{-n\alpha
q}}{K^m\tau^{mn}}.
\end{eqnarray*}

Observing that for each $n$, the total number of the events
$\cA_{z,\tau}^{n,e}$  ($z,z+e\in\cS_n$) is not greater than
$2^{c(d)n}$. Hence the probability of the union
$\cA_\tau^n=\cup_{z,z+e\in
S_n}(\cup_{\|e\|_\infty=1}\cA_{z,\tau}^{n,e})$ meets the estimate
\begin{eqnarray*}
{\mathbb P}(\cA_\tau^n)\leq C(m,T)
\frac{2^{-nm\alpha}}{K^m\tau^{mn}}2^{c(d)n}\leq C(m,T)K^{-m}
\Big(\frac{2^{c(d)}}{(2^\alpha\tau)^m}\Big)^n.
\end{eqnarray*}

One chooses $\tau=2^{-\alpha/2}$, $m=3c(d)/\alpha\vee 1$, then the
probability of the event $\cA=\cup_{n\geq1}\cA_\tau^n$ can  be
calculated that
\begin{eqnarray}\label{2.27}
{\mathbb P}(\cA)\leq C(T,d)K^{-m}.
\end{eqnarray}

For every point $x\in [0,1]^d$, we have $x=\sum_{i=0}^\infty
e_i2^{-i}$ ($\|e_i\|_\infty\leq 1$). Denote
$x_k=\sum_{i=0}^ke_i2^{-i}$. For any $\omega\overline{\in}\cA$,
we have $|\xi_t(x_{k+1})-\xi_t(x_k)|< K\tau^{k+1}$, which suggests
that
\begin{eqnarray}\label{2.28}
|\xi_t(x)-\xi_t(x_0)|\leq \sum_{k=0}^\infty|\xi_t(x_{k+1})-\xi_t(x_k)|<
K\sum_{k=0}^\infty\tau^{k+1}\leq CK,
\end{eqnarray}
where we have fetched $\tau=2^{-\alpha/2}$.

Set $\xi_1=\sup_{(t,x)\in [0,T]\times [0,1]^d}|\xi_t(x)-\xi_t(x_0)|$, then
for any $0<\gamma<m$, it yields that
\begin{eqnarray}\label{2.29}
{\mathbb E} |\xi_1|^\gamma=\gamma\int_0^\infty  \lambda^{\gamma-1}{\mathbb P}(\xi_1\geq \lambda)d\lambda
=\gamma\int_0^{CK}\lambda^{\gamma-1}{\mathbb P}(\xi_1\geq \lambda)d\lambda+\gamma\int_{CK}^\infty
\lambda^{\gamma-1}{\mathbb P}(v\geq \lambda)dr.
\end{eqnarray}
According to (\ref{2.28}) and (\ref{2.27}), from (\ref{2.29}) one
finishes at
\begin{eqnarray*}
{\mathbb E} |\xi_1|^\gamma\leq (CK)^\gamma+C(T,d) \gamma\int_{CK}^\infty
\lambda^{\gamma-1-m}d\lambda \leq (CK)^\gamma +C(T,d)\gamma K^{\gamma-m},
\end{eqnarray*}
which hints that
\begin{eqnarray}\label{2.30}
{\mathbb E} \sup_{(t,x)\in [0,T]\times [0,1]^d}|\xi_t(x)|^\gamma \leq C(\gamma)\Big[{\mathbb E} |\xi_1|^\gamma +{\mathbb E} \sup_{0\leq t\leq T}|\xi_t(x_0)|^\gamma\Big]\leq
C(T,d),
\end{eqnarray}
where in the last inequality we have used the estimate (\ref{2.22}).

Hence, we conclude that
\begin{eqnarray*}
{\mathbb P}(\|\xi_\cdot(\cdot)\|_{L^\infty((0,T)\times [0,1]^d)}<\infty)=1,
\end{eqnarray*}
and by Step 1 thus
\begin{eqnarray}\label{2.31}
{\mathbb P}(\|\xi_\cdot(\cdot)\|_{L^\infty((0,T)\times B_R)}<\infty)=1,
\quad \forall \ R>0.
\end{eqnarray}

Therefore
\begin{eqnarray*}
&& \quad  \|\nabla_x(u_0(X^{-1}))\|_{L^\infty((0,T)\times B_R)}\leq C(\omega)
\|\nabla_xu_0\|_{L^\infty(B_{C(\omega)})}<\infty,  \ \ {\mathbb P}-a.s., \ \ when \ \ r=\infty, \cr\cr&&
\int_{B_R}|\nabla_x(u_0(X^{-1}(t,x)))|^rdx\leq C(\omega)
\int_{B_{C(\omega)}}|\nabla_xu_0(x)|^rdx<\infty, \ \ {\mathbb P}-a.s. , \ \ when \ \ r\in [1,\infty).
\end{eqnarray*}
From above estimated, one accomplishes the proof for Case 1.

\vskip2mm\par
\textbf{$\bullet$ Case 2:} $p\in [2,\infty), q\in (2,\infty)$, $\frac{2}{q}+\frac{d}{p}<1$.
\vskip2mm\par
The proof is the same as the proof for Case 1 since now for every $t\in[0,T]$, $\gamma(t,x)=x+U(t,x)$ forms a non-singular diffeomorphism of $\mathcal{C}^2$. Moreover, $\gamma$ and $\gamma^{-1}$ have a bounded first and second spatial derives and $\nabla^2 \gamma$ is globally H\"{o}lder continuous.

(ii) To reach our aim, we firstly notice the following two facts:
$m$ is arbitrary in (\ref{2.26}), and (\ref{2.30}) holds true for every $0<\gamma<m$. Then by the statements in \textbf{The second proof}, \textbf{Step 1}, for every $a\geq 1$, we derive that
\begin{eqnarray}\label{2.32}
{\mathbb E}|\nabla X^{-1}(\cdot,\cdot)|_{L^\infty((0,T)\times B_R)}^a \leq C(T,d,R,a), \quad \forall \ \ R>0.
\end{eqnarray}
If $|\nabla u_0|\in L^\infty({\mathbb R}^d)$, from (\ref{2.32}), for every $a\geq 1$, and every $R>0$, one arrives at
\begin{eqnarray*}
{\mathbb E}\sup_{(t,x)\in [0,T]\times B_R}\|\nabla_xu(t,x)\|^a&=&{\mathbb E} \sup_{(t,x)\in [0,T]\times B_R} \|\nabla_xu_0(X^{-1}(t,x))\nabla X^{-1}(t,x)\|^a\cr\cr&\leq &C(T,d,R,a) \|\nabla_xu_0\|_{L^\infty({\mathbb R}^d)}^a.
\end{eqnarray*}
So we finish the proof. $\Box$

\vskip2mm\noindent
\textbf{Remark 2.2.} (i) From '\textbf{The first proof}', one also achieves that: let $\beta\in (0,1)$ that $u_0\in \mathcal{C}_b^\beta({\mathbb R}^d)$, $b\in L^\infty(0,T;\mathcal{C}_b^\beta({\mathbb R}^d;{\mathbb R}^d))$, $\div b\in L^\infty((0,T)\times {\mathbb R}^d)$, there exists a unique bounded random field $u$ which takes value in  $L^\infty(0,T;\mathcal{C}_b^\beta({\mathbb R}^d)) \ \ {\mathbb P}-a.s.$ such that for every $\varphi\in\mathcal{C}_0^\infty({\mathbb R}^d)$, $\int_{{\mathbb R}^d}\varphi(x)u(t,x)dx$ has a continuous modification which is an $\mathcal{F}_t$-semimartingale and for every $t\in [0,T]$, (\ref{1.2}) holds. Moreover, $u(t,x)=u_0(X^{-1}(t,x))$. On the other hand, if one applies '\textbf{The second proof}' from Step 1 to Step 4, it yields that for every $R>0$, $u\in \cap_{a\geq 1}L^a(\Omega;L^\infty(0,T;\mathcal{C}_b^\beta(B_R)))$.

(ii) For irregular drift coefficient, the existence and uniqueness of stochastic weak solutions of (\ref{1.1}) can be seen in \cite{CO}. The existence and non-existence of $L^\infty\cap W^{1,r}$ solutions for the deterministic equations of (\ref{1.1}) can be found in \cite{CLR2,Jab}.

 \vskip2mm\noindent
\textbf{Remark 2.3.}  (i)  The estimate
(\ref{2.27}) for the tail probability is inspired by \cite{KNP}. In
\cite{KNP}, Kuksin, Nadirashvili and Piatnitski argued
\begin{eqnarray*}
du(t,x)=Au(t,x)dt+f(t,x)dB(t), \ t>0, \ x\in Q, \ \ u(t,x)|_{t=0}=0,
\end{eqnarray*}
where $Q$ is a bounded domain,  by estimating the tail probability,
they gained a space and time H\"{o}lder estimates for solutions,
i.e. ${\mathbb E} \|u\|_{\mathcal{C}^\theta(Q_T}^\gamma <\infty$. Here, we get an analogue of \cite{KNP}. For more details, one also refers to \cite{Kuk}.

(ii) For simplicity, here we only discuss the noise given by
$\sum_{i=1}^d
\partial_{x_i}u\circ\dot{B}_i(t)$, however, the method is appropriate for the stochastic transport equation
below
\begin{eqnarray*}
\partial_tu(t,x)+b(t,x)\cdot\nabla_x u(t,x)
+\sum_{i=1}^d\sum_{j=1}^n\partial_{x_i}u(t,x)G_{i,j}(t,x)\circ\dot{B}_j(t)=0,
\ (\omega,t,x)\in \Omega\times(0,T)\times{\mathbb R}^d,
\end{eqnarray*}
here $1\leq i\leq  d, \ 1\leq j\leq n$, $d,n\in\mN$. But now one should replace the SDE (\ref{2.1}) by
\begin{eqnarray*}
dX(t)=b(t,X(t))dt+G(t,X(t))dB(t),  \ \
 X(0)=x\in {\mathbb R}^d,
 \end{eqnarray*}
where $G(t,x)=(G_{i,j}(t,x))\in {\mathbb R}^{d\times n}$, $B=(B_1,B_2,_{\cdots}, B_n)^\top$ is a standard $n$-dimensional Brownian motion.  For more details about above SDE, the authors can consult to \cite{Duan}.

(iii) How to extend the present result to the nonlinear scalar conservation law
\begin{eqnarray}\label{2.33}
\left\{
  \begin{array}{ll}
\partial_tu(t,x)+\div F(u(t,x))=0, \ (t,x)\in (0,T)\times {\mathbb R}^d, \\
u(t,x)|_{t=0}=u_0(x), \  x\in{\mathbb R}^d,
  \end{array}
\right.
\end{eqnarray}
may be an very interesting problem. Observing that when $F$ and $u$ are smooth, we rewrite (\ref{2.33}) by
\begin{eqnarray*}
\left\{
  \begin{array}{ll}
\partial_tu(t,x)+F^\prime(u)\cdot\nabla u(t,x)=0, \ (t,x)\in (0,T)\times {\mathbb R}^d, \\
u(t,x)|_{t=0}=u_0(x), \  x\in{\mathbb R}^d.
  \end{array}
\right.
\end{eqnarray*}
With this formulation in mind, it is desirable that the coefficients are random and have no smoother than the solutions (at least as smooth as solutions).
As stated in \cite{FGP1}, it would be a difficult problem.
However, for linear transport equation, when $b$ is only depends on the random perturbation $B$, there has some celebrate works such as see \cite{DR}.

\section{Proof of Theorem \ref{the1.2}}\label{sec4}
\setcounter{equation}{0}
\vskip2mm\noindent
\textbf{Proof.} The existence, uniqueness as well as obvious representation of weak solutions can be seen in
DiPerna and Lions \cite{DL}, we omit some details. Now let us show the non-existence of strong solutions.

For simplicity, we assume $b$ is time independent and $d=2$. Now we rewrite $x$ by $(x,y)\in{\mathbb R}^2$. Let $b_1(x)$
and $b_2(y)$ be defined as the following:
\begin{eqnarray*}
b_1(x)=\left\{ \begin{array}{ll}
x^{\frac{3}{4}},  \quad 0<x<1 \\ x^{-\frac{3}{4}}, \quad x\geq 1, \\ \ \ 0, \quad otherwise,\end{array}
\right.
b_2(y)=\left\{ \begin{array}{ll}
\frac{y}{1+y^2},  \quad y\geq 0, \\ \ \
 0, \quad otherwise.\end{array}
\right.
\end{eqnarray*}
Then $b_1,b_2\in L^3({\mathbb R})$ and $0\leq b_1,b_2\leq 1$. Notice that
\begin{eqnarray*}
\frac{d}{dx}b_1(x)=\left\{ \begin{array}{ll}
\frac{3}{4}x^{-\frac{1}{4}},  \quad 0<x<1 \\ -\frac{3}{4}x^{-\frac{7}{4}}, \quad x> 1, \\ \ \ 0,
\quad otherwise,\end{array}
\right.
\frac{d}{dy}b_2(y)=\left\{ \begin{array}{ll}
\frac{1-y^2}{(1+y^2)^2},  \quad y\geq 0, \\ \ \
 0, \quad otherwise,\end{array}
\right.
\end{eqnarray*}
we conclude that $b_1,b_2\in W^{1,3}({\mathbb R})$ and $-1/8\leq b_2^\prime\leq 1$, but $\sup b_1^\prime=\infty$.

We define $b(x,y)=(0,b_1(x)b_2(y))$, then $b\in W^{1,3}({\mathbb R}^2)$,
$\div b(x,y)=b_1(x)b_2^\prime(y)$ and $\div b\in [-1/8,1]$. Consider the ODE below
\begin{eqnarray*}
\frac{d}{dt}X(t)=0, \ \frac{d}{dt}Y(t)=b_1(X(t))b_2(Y(t)),
\ X(0)=x\geq 0, \ Y(0)=y\geq 0,
\end{eqnarray*}
we gain
\begin{eqnarray}\label{3.1}
X(t,x)=x, \ Y(t,x,y)=g^{-1}(g(y)e^{2b_1(x)t}),
\end{eqnarray}
where $g(y)=e^{y^2}y^2 \ (y\geq 0)$, $g^{-1}$ is the inverse of $g$.

From (\ref{3.1}), then
\begin{eqnarray*}
\frac{\partial(X,Y)}{\partial(x,y)}=
\left(  \begin{array}{cc} 1 & (g^{-1})^\prime(g(y)e^{2b_1(x)t})g(y)e^{2b_1(x)t}2b_1^\prime(x)t \\                                        0 & (g^{-1})^\prime(g(y)e^{2b_1(x)t})g^\prime(y)e^{2b_1(x)t} \\
\end{array}  \right),
\end{eqnarray*}
and for every real number $R>0$,
\begin{eqnarray}\label{3.2}
(X(t),Y(t))([0,R]\times [0,R])\supset [0,R]\times [0,R].
\end{eqnarray}

With the inverse function theorem, it yields that
\begin{eqnarray}\label{3.3}
\Big(\frac{\partial(X(t),Y(t))}{\partial(x,y)}\Big)^{-1}(x,y)=
\left(  \begin{array}{cc} 1 & -(g^\prime(y))^{-1}g(y)2b_1^\prime(x)t \\
 0 & (g^\prime(g(y)e^{2b_1(x)t})(g^\prime(y))^{-1}e^{-2b_1(x)t} \\
\end{array}  \right).
\end{eqnarray}
Combining the fact $u(t,x,y)=u_0((X,Y)^{-1}(t,x,y))$ and (\ref{3.2}), for every $t>0$, every $R>0$, one ends up with
\begin{eqnarray}\label{3.4}
&&\int_{[-R,R]\times [-R,R]}|\nabla_{x,y}(u_0((X,Y)^{-1}(t,x,y)))|^3dxdy\cr\cr&\geq&
\int_{[0,R]\times [0,R]}|\nabla_{x,y}(u_0((X,Y)^{-1}(t,x,y)))|^3dxdy\cr\cr&\geq&
\int_{[0,R]\times [0,R]}\!\!|\nabla_{x,y}u_0(x,y)|^3
\Big\|\Big(\frac{\partial(X,Y)}{\partial(x,y)}\Big)^{-1}\Big
\|^3\!\!\exp(\int^t_0\!\!\div
b(X(r,x),Y(r,x,y))dr)dxdy.
\end{eqnarray}

With the aid of Louville's theorem, (\ref{3.3}) and noticing that $-1/8\leq \div b\leq 1$,  from (\ref{3.4}), one arrives at
\begin{eqnarray}\label{3.5}
&&\int_{[-R,R]\times [-R,R]}|\nabla_{x,y}(u_0((X,Y)^{-1}(t,x,y)))|^3dxdy\cr\cr&\geq&
\exp(-\frac{t}{8})
\int_{[0,R]\times [0,R]}|\nabla_{x,y}u_0(x,y)|^3|(g^\prime(y))^{-1}g(y)2b_1^\prime(x)t|^3
dxdy\cr\cr&\geq&\exp(-\frac{t}{8})t^3
\int_{[0,R]\times [0,R]}|\nabla_{x,y}u_0(x,y)|^3\Big(\frac{y}{1+y}\Big)^3|b_1^\prime(x)|^3dxdy,
\end{eqnarray}
where in the last inequality we have used
$$
(g^\prime(y))^{-1}g(y)=\frac{y}{2(1+y)}, \quad \forall \ y\geq 0.
$$
If one chooses $u_0(x,y)=u_{0,1}(x)u_{0,2}(y)$ such that $u_{0,1},u_{0,2}\in W^{1,3}({\mathbb R})$ and
$u_{0,1}^\prime(x)\approx x^{-1/4}$ near $0+$, then
\begin{eqnarray*}
&&\int_{[-R,R]\times [-R,R]}|\nabla_{x,y}(u_0((X,Y)^{-1}(t,x,y)))|^3dxdy
\cr\cr&\geq&\exp(-\frac{t}{8})t^3\int_0^R|u_{0,2}^\prime(y)|^3\Big(\frac{y}{1+y}\Big)^3dy
\int_0^R |u_{0,1}^\prime(x)|^3|b_1^\prime(x)|^3dx
\cr\cr&\geq&C\exp(-\frac{t}{8})t^3
\int_0^Rx^{-\frac{3}{4}}|b_1^\prime(x)|^3dx
\cr\cr&\geq&C\exp(-\frac{t}{8})t^3
\int_0^{\epsilon}x^{-\frac{3}{4}}|b_1^\prime(x)|^3dx =\infty,
\end{eqnarray*}
where $\epsilon>0$ is a small enough positive real number. From this we complete the proof. $\Box$

\vskip2mm\noindent
\textbf{Remark 3.1.} Based upon Theorem \ref{the1.1} and Theorem \ref{the1.2}, we deduce that noise can prevent singularities and here let us give an explanation. From above construction we know that if one chooses $u_0(x,y)=u_{0,1}(x)u_{0,2}(y)$ such that $u_{0,1},u_{0,2}\in W^{1,3}({\mathbb R})$ and
$u_{0,1}^\prime(x)\approx x^{-1/4}$ near $0+$, for the deterministic equation, then
\begin{eqnarray}\label{3.6}
\int_{[-R,R]\times [-R,R]}|\nabla_{x,y}u(t,x,y)|^3dxdy\geq C(t)
\int_0^Rx^{-\frac{1}{4}}|b_1^\prime(x)|^3dx =\infty.
\end{eqnarray}
However, when one deals with the stochastic equation, the characteristic equation becomes into
\begin{eqnarray*}
dX(t)=dB_1(t), \ dY(t)=b_1(X(t))b_2(Y(t))dt+dB_2(t),
\ X(0)=x, \ Y(0)=y.
\end{eqnarray*}
Similar calculations implies the estimate
\begin{eqnarray}\label{3.7}
&&\int_{[-R,R]\times [-R,R]}|\nabla_{x,y}u(t,x,y)|^3dxdy\cr\cr&\leq& C(t)
\Big[\int_0^Rx^{-\frac{3}{4}}|b_1^\prime(x+B_1(t))|^3dx+1\Big], \quad {\mathbb P}-a.s. \cr\cr&\leq& C(t)
\Big[\int_0^1x^{-\frac{3}{4}}|x+B_1(t))|^{-\frac{3}{4}}dx+1\Big],\quad {\mathbb P}-a.s.
\end{eqnarray}
for every $R>0$.

Since for every $t>0$, $B_1(t)$ has a normal distribution with expected value $0$ and variance $t$,
\begin{eqnarray}\label{3.8}
{\mathbb E}\int_0^1x^{-\frac{3}{4}}|x+B_1(t))|^{-\frac{3}{4}}dx=\int_0^1x^{-\frac{3}{4}}dx\int_{{\mathbb R}}
|y|^{-\frac{3}{4}}\frac{1}{\sqrt{2\pi t}}e^{-\frac{(y-x)^2}{2t}}dy.
\end{eqnarray}
In view of Cauchy-Schwarz's inequality,
\begin{eqnarray*}
2|xy|\leq \frac{y^2}{2}+2x^2,
\end{eqnarray*}
then
\begin{eqnarray}\label{3.9}
-(y-x)^2\leq -y^2-x^2+\frac{y^2}{2}+2x^2=-\frac{y^2}{2}+x^2.
\end{eqnarray}
Combining (\ref{3.8}) and (\ref{3.9}), we reach at
\begin{eqnarray*}
{\mathbb E}\int_0^1x^{-\frac{3}{4}}|x+B_1(t))|^{-\frac{3}{4}}dx
\leq \int_0^1x^{-\frac{3}{4}} e^{\frac{x^2}{2t}} dx
\int_{{\mathbb R}}|y|^{-\frac{3}{4}}\frac{1}{\sqrt{2\pi t}}e^{-\frac{y^2}{4t}}dy<\infty,
\end{eqnarray*}
which hints that
\begin{eqnarray*}
\int_0^1x^{-\frac{3}{4}}|x+B_1(t))|^{-\frac{3}{4}}dx<\infty,\quad {\mathbb P}-a.s..
\end{eqnarray*}
From (\ref{3.7}), therefore
\begin{eqnarray*}
\int_{[-R,R]\times [-R,R]}|\nabla_{x,y}u(t,x,y)|^3dxdy<\infty,\quad {\mathbb P}-a.s..
\end{eqnarray*}

\vskip2mm\noindent
\textbf{Remark 3.2.} (i) When $W^{1,p}$ is replaced by a complete topological vector space $S$, such that
\begin{eqnarray*}
BV_{loc}\cap L^\inf({\mathbb R}_t\times{\mathbb R}^d_x)\subset S({\mathbb R}_t\times{\mathbb R}^d_x)\subset\subset L^1_{loc}\cap
L^\inf({\mathbb R}_t\times{\mathbb R}^d_x).
\end{eqnarray*}
The non-existence of (weak) solutions for (\ref{1.6}) has been
proved by Crippa and De Lellis \cite{CD} for $d\geq 3$ (similar
question can be seen in \cite{Bre,Dep}). But, when one deals with
(\ref{1.1}), one may give a positive answer for existence and
uniqueness. We leave this topic in a further work.

(ii) For linear transport equation (\ref{1.6}), if one replaces $L^\infty(0,T;W^{1,p}({\mathbb R}^d))$
by $H^{1/2}\cap L^\infty((0,T)\times{\mathbb R}^d)$, then using Leibniz's rule, Colombini, Luo and Rauch \cite{CLR1}
obtained the uniqueness of $H^{1/2}\cap L^\infty((0,T)\times{\mathbb R}^d)$ solutions. Without Sobolev regularity,
the authors gave a counter example on uniqueness. Thus there is some proper space $S$,
$H^{1/2}\cap L^\infty((0,T)\times{\mathbb R}^d)\subset S\subset L^\infty((0,T)\times{\mathbb R}^d)$, such that the Cauchy
problem (\ref{1.4}) is well-posed in $S$. Moreover, using stochastic regularization, we may gain a bigger
space $S_1$, so that in $S_1$ the stochastic equation is well-posed.

\end{document}